\documentstyle[12pt]{article}

\advance\hoffset by -7mm
\makeatletter
\@addtoreset{equation}{section}
\makeatother

\input amssym.def  
\input amssym.tex

\renewcommand{\Pi}{{\Bbb T}}

\renewcommand{\title}[1]{\null\vspace{25mm}

\noindent{\Large{\bf #1}}\vspace{10mm}

\noindent {\large By }}
\newcommand{\authors}[1]{\noindent{\large #1}\vspace{3mm}

}
\newcommand{\address}[1]{\noindent #1\vspace{5mm}

}
\renewcommand{\abstract}[1]{\vspace{19mm}

\noindent{\small{\em Abstract.} #1}\vspace{2mm}

}

\begin{document}

\title{Symplectic twist maps without conjugate points}
\authors{M L Bialy}
\address{School of Mathematical Sciences, Tel-Aviv University,
69978 Tel-Aviv,  Israel\\
Email:  bialy@math.tau.ac.il}
\authors{and R S MacKay}
\address{Mathematics Institute, University of
Warwick, Coventry CV4 7AL, UK\\
Email:  mackay@maths.warwick.ac.uk}
\date{version of \today}

\abstract{
For sequences of symplectic twist maps without conjugate points, an
invariant Lagrangian subbundle is constructed.  This allows one to deduce
that absence of conjugate points is a rare property in some classes of map.
}

\bigskip
\noindent{This work was partially supported by EPSRC grant
GR/M11349.}

\section{Introduction and results}

In this paper we construct an analogue of L Green's invariant
subbundles for the case of discrete variational principles related
to the dynamics of sequence of symplectic twist maps of $T^*
\Pi^d$. Such a construction was first performed by L Green
\cite{[G]} for Riemannian geodesic flows but has turned out to be much
more general. For example, it can be extended to optical
Hamiltonian flows \cite{[C-I]}.  The construction of
invariant subbundles is very useful in many examples of the
so-called Hopf-type rigidity.

In particular, we apply L Green's construction to the so-called
Frenkel-Kontorova variational problem which is related to a
sequence of generalized standard maps.  We prove a result which
can be seen as an analogue of a rigidity result of Knauf and
Croke-Fathi which was proved for conformally flat Riemannian
metrics \cite{[K],[C-F]}.

In the discrete time case  Hopf rigidity was established first for
convex plane billiards \cite{[B1],[W]}.

There are still very many problems related to the rigidity and
integrability of twist maps and we hope that our results will be
useful for their solutions.

Let us introduce the setting (see also the recent book by Chris
Gole \cite{[Go]} for a detailed exposition).

For each $n \in \Bbb Z$, let $S_n: \Bbb R^d \times \Bbb R^d
\rightarrow \Bbb R$ be a $C^2$-smooth function satisfying the
following:
\begin{eqnarray}
1. && S_n \quad {\rm is} \quad \Bbb Z^d -{\rm periodic:}\quad S_n (q+e, Q+e)
= S_n(q,Q) \nonumber\\
&& {\rm for \; any}\; (q,Q) \in \Bbb R^d \times \Bbb R^d \; {\rm and}
\; e \in \Bbb Z^d.
\label{11}
\end{eqnarray}

2. $S_n \quad {\rm satisfies\; the\; uniform\; twist\; condition:}
{\rm \; for \; any\;} \xi \in \Bbb R^d {\rm \;the\; quadratic\;
form}$
\begin{equation}  \sum_{i, j} {\partial^2 S_n (q,Q) \over \partial q_i
\partial Q_j} \xi_i \xi_j \leq - K ||\xi||^2 {\rm \; for\; a\;
positive\; constant\;} K. \label{12}
\end{equation}
Such a function defines two closely related objects.

The first is the variational functional defined on the sequences
$\left\{ q_n \right\}, n \in \Bbb Z$,
\begin{equation}
F\left( \left\{ q_n \right\} \right) =  \sum^{+ \infty}_{n= -
\infty} S_n \left( q_n, q_{n+1} \right) . \label{13}
\end{equation}
The functional is a formal sum but the extremals are well defined
and satisfy the equations
\begin{equation}
\partial_2 S_{n-1} \left(q_{n-1}, q_n \right) + \partial_1 S_n
\left( q_n, q_{n+1} \right) = 0 \; {\rm for\; all\;} n \in \Bbb Z.
\label{14}
\end{equation}

The second object is the symplectic diffeomorphism $T_n$ of $T^*
\Pi^d$ generated by the function $S_n$.  In the standard
coordinates $(p,q)$ it is given by the following implicit formula
\begin{equation}
T_n(p,q) = (P, Q) \; {\rm if\;} P = + \partial_2 S_n (q,Q), p = -
\partial_1 S_n (q,Q). \label{15}
\end{equation}
Here and throughout the paper $\partial_1, \partial_2$ stand for the
derivatives with respect to the $q_i, Q_j$ variables respectively.

We refer the reader to \cite{[Go]}, \cite{[He]} and \cite{[M-M-S]} for
general theory of symplectic twist maps --- note
that in eq.(1.2) we follow \cite{[M-M-S]}'s choice of twist condition
rather than either of those of \cite{[He]}.

The basic example for us will be

\noindent{\em Example 1}  Let $S_n = {1 \over 2} || Q-q||^2 +
V_n(q)$ where $V_n$ is a $\Bbb Z^d$-periodic smooth function (called
the potential).  In
this case we shall call $F$ a Frenkel-Kontorova functional. In what
follows we will assume that the sequence of the potential
functions $V_n$ depends either periodically on $n$, or $V_n$
vanishes for all but finitely many values of $n$. The
corresponding map $T_n$ is a generalized standard map of $T^*
\Pi^d$:
$$
T_n: (p,q) \mapsto (p + \nabla V_n (q), p+q+\nabla V_n(q)).
$$
It is important to notice that in this case for any $n$, $T_n$ can be
considered as acting on $\Pi^{2d}$ and not just on $T^*\Pi^d$; this
follows from the fact that for any $e \in \Bbb Z^d$,
$$
T_n(p+e, q) = ( P+e, Q+e) .
$$

The correspondence between the extremals of the functional $F$ and
the orbits of the sequence $T_n$ is the following. Let a sequence
$\left\{ q_n \right\}$ be an extremal for $F$. Let $p_n = -
\partial_1 S(q_n, q_{n+1})$ and form the sequence $\{x_n = (p_n,
q_n)\}$. Then $\{x_n\}$ is an orbit of the evolution, i.e. $
T_n(x_n)=x_{n+1}$. Conversely, if $\{x_n = (p_n, q_n)\}$ is an
orbit then the corresponding sequence $\{q_n\}$ is extremal for
the variational principle written above.

Similarly, invariant fields along the orbits of $\{T_n\}$
correspond to the so-called Jacobi fields along the extremals. For
an orbit $\{x_n\}$, let $\zeta_n \in T_{x_n} T^* \Pi^d$ be a
tangent vector at $x_n = (p_n, q_n)$; then the field $\{\zeta_n\}$
is invariant under the derivative
$T_*$, i.e. $(T_n)_*(\zeta_n)=\zeta_{n+1}$, if
and only if the vectors $\xi_n = \pi_* (\zeta_n)$ satisfy the
Jacobi equation (here $\pi: (p,q) \mapsto q$ is the canonical
projection):
\begin{equation}
b^T_{n-1} \xi_{n-1} + a_n \xi_n + b_n \xi_{n+1} = 0
\label{16}
\end{equation}
with the matrices $$b_n = \partial_{12} S_n(q_n, q_{n+1}), a_n =
\partial_{11} S_n(q_n, q_{n+1}) + \partial_{22} S_{n-1} (q_{n-1},
q_n)$$ (the symbols $\partial_{11}S, \partial_{12}S,
\partial_{22}S$ denote the matrices of second derivatives of $S$).

We will use the following definition first introduced for the discrete
case in \cite {[B1]}.

\noindent{\em Definition}.  Two points of the extremal configuration
$\left\{ q_n \right\}$ are called {\em conjugate} if there exists a
non-trivial Jacobi field $\xi_n$ vanishing at these two points.

Denote by $R^n_m$ the evolution transformation, i.e.

$ R^n_m=T_{n-1}\circ\ldots\circ T_m$, for ${n>m}$, $R^m_m=Id$ and
$R^n_m=(R^m_n)^{-1}$, for $n<m$.

With the above correspondence one can interpret the definition geometrically
by saying that $q_m$ and $q_n$, for $m<n$, are conjugate if
$$
(R^{n-1}_m)_* \left({\cal V} (x_m) \right) \bigcap {\cal V} (x_n)
\not= \left\{ 0 \right\}
$$
where ${\cal V}(x)$ denotes the vertical subspace at $x$ and $x_n
=R_0^n( x_0)$ is the orbit corresponding to $\left\{ q_n
\right\}$.

\newtheorem{theorem}{Theorem}
\begin{theorem}
\label{T1} If  none of the extremals of the functional $F$ have
conjugate points then for every {n} there exists a field ${\cal
W}_n$ of Lagrangian subspaces ${\cal W}_n(x) \subseteq T_x T^*
\Pi^d$ depending measurably on $x$ and such that
\begin{enumerate}
\item Invariance: $(T_n)_* {\cal W}_n(x) ={\cal W}_{n+1}(T_nx)$
\item At every point $x, {\cal W}_n(x)$ is transversal to the
vertical subspace
  ${\cal V}(x)$.
\end{enumerate}
\end{theorem}
We shall use a partial order $\leq$ on the subset of Lagrangian
subspaces which are transversal to the vertical one, defined as
follows. To every such subspace ${\cal L}(x)$ corresponds a
symmetric matrix $L$, by ${\cal L}(x) = \left\{ \xi: dp(\xi)=L\ 
dq(\xi) \right\}$. Given two subspaces ${\cal L}_{1,2}$ we say
${\cal L}_1 \leq {\cal L}_2$ if $L_1 \leq L_2$, i.e.~$L_2 - L_1$
is non-negative.

\begin{theorem}
\label{T2} If none of the extremals of the functional $F$ have
conjugate points then for the fields ${\cal W}_n(x)$ the following
holds
\begin{enumerate}
\item $(T_{n+1}^{-1})_* {\cal V}(T_nx) \le {\cal W}_n(x) \le
(T_{n-1})_* \left({\cal V}(T_n^{-1}x) \right)$,

or in terms of the matrices this reads

$-\partial_{11}S_n (q,q_{+}) \leq W_n(x) \leq
  \partial_{22} S_{n-1}(q_{\_}, q)$,

for all $x$ where $q=\pi(x), q_{\_}=\pi(T^{-1}_{n-1}x),
q_{+}=\pi(T_n x)$.

\item The following inequality holds true for all $x$

$W_{n+1} (T_n x) - W_n (x) \leq
\partial_{11} S_n (q, q_{+}) +
\partial_{22}
  S_n (q, q_{+}) + \partial_{12} S_n (q, q_{+}) + \partial_{21} S_n (q,
  q_{+})$

  with equality in only the case when

  $\partial_{12} S_n(q, q_{+}) =
  \partial_{21}S_n (q, q_{+})$ and $W_{n+1} (T_nx) = \partial_{22}S_n (q, q_{+})
  + \partial_{21}S_n (q, q_{+})$

  and $W_n(x) = - \partial_{11} S_n (q,
  q_{+}) - \partial_{12} S_n (q, q_{+})$.
\end{enumerate}
\end{theorem}

As an application of this to Frenkel-Kontorova functionals we
obtain

\begin{theorem}
\label{T3} Consider the Frenkel-Kontorova functional with a
sequence of potential functions $V_n$ which
is either periodic in $n$ or has all but finitely many of the $V_n$s
constant functions.  Then either there exist extremals with
conjugate points or all the potential functions are constants.
\end{theorem}

The next section contains necessary preliminaries about Jacobi fields
in the discrete case.  We prove the theorems in section 3.  Discussion
and open questions conclude the paper.

\noindent{\large \bf Acknowledgements}

\noindent This paper was started while the first author was
visiting the second one at the Nonlinear Centre, University of
Cambridge in 98/99 in the framework of a joint project EPSRC grant.
We are grateful to the EPSRC for their support and to the University of
Cambridge for their hospitality.

\section{Nonsingular Jacobi fields}

In this section we prove first that the assumption that no
extremal has conjugate points implies that each extremal is in
fact a strict local minimum configuration. As a consequence of this we
construct a special non-singular solution of the {\em matrix}
Jacobi equation. The first fact is stated as
\newtheorem{lemma}{Lemma}
\begin{lemma}
\label{L1} If all the extremals of $F$ have no conjugate points
then each is a strict local minimum between any two of its points.
\end{lemma}

\noindent{\bf Proof of Lemma {\ref{L1}}}

Let $\left\{ q_n\right\}, n \in \Bbb Z$, be an extremal.  For $M \le N$,
denote
$$
F_{MN} (u_M, \ldots, u_N) = S_{M-1}(q_{M-1}, u_M) +
\sum^{N-1}_{n=M} S_n(u_n, u_{n+1}) + S_N(u_N, q_{N+1}).
$$
We claim that the matrix $\delta^2 F_{M,N}$ of second variation of
$F_{MN}$ is positive definite.  To prove this, note that by a
simple calculation it has the following block matrix form:
\begin{eqnarray}
\left (\begin{array}{cccc}
a_M&b_M& &0\\
b_{M}^T&a_{M+1}&\ddots&\\
&\ddots&\ddots&b_{N-1}\\
0&&b_{N-1}^T&a_N
\end{array}\right)
\end{eqnarray}
with the matrices $a_i, b_i$ introduced in eq.(1.6).  It follows
that the kernel of this matrix consists exactly of the Jacobi
fields vanishing at $q_{M-1}$ and $q_{N+1}$.  Thus by the non
conjugacy assumption, the matrix is non-degenerate.  But then it
has to be positive definite by the fact that it depends
continuously on the configuration (and so its signature is
constant) and there always exist segments which minimize the
functional (a consequence of (1.1),(1.2)) and so have positive
definite second variation (see for example \cite{[Go]} for the
proof). This completes the proof of the lemma. $\Box$

Note that as a consequence, every orbit is a global minimum between any two
of its points, though we do not need this fact.

Let us consider a minimal configuration $\left\{ q_n \right\}, n \in
\Bbb Z$.  For given $k \in \Bbb Z$,
define a {\em matrix} solution of the Jacobi eq.(1.6)
$\xi^{(k)}_{n}$ such that $\xi^{(k)}_{k}=0$ and $\xi^{(k)}_{k+1}$ is
invertible, by iteration from this pair.
Then by the no conjugate points assumption, all $\xi^{(k)}_{n}$  are
invertible $(n \not= k)$ and hence
\begin{equation}
A^{(k)}_n = - b_n \xi^{(k)}_{n+1} \left[ \xi^{(k)}_n \right]^{-1} (n>k)
\label{22}
\end{equation}
are defined and do not depend on the choice of
$\xi^{(k)}_{k+1}$.  Moreover one can easily see that
\begin{equation}
A^{(k)}_{k+1} = a_{k+1}, \mbox{ and  for } n>k, \quad
A^{(k)}_{n+1} = a_{n+1} -b^T_n \left[ A^{(k)}_n \right]^{-1} b_n .
\label{23}
\end{equation}
In particular all the $A^{(k)}_n$ are symmetric.  A crucial
observation for us is that all these matrices are in fact positive
definite.  Indeed, if on the contrary, for some $m>k$, $A^{(k)}_m$
is not positive definite then for some vector $\eta \ne 0$,
$<A^{(k)}_m \eta, \eta >\ \leq 0$.  Then define the segment of
Jacobi field
$$
\eta_n = \xi^{(k)}_n \left[ \xi^{(k)}_m \right]^{-1} \eta, \quad k
\leq n \leq m ;
$$
for $n=k$ and $n=m$ we have $\eta_k = 0$, $\eta_m = \eta$.
One can easily compute the value of the quadratic form $\delta^2
F_{k+1, m}$ on the variation $(\eta_{k+1}, \dots, \eta_m)$.
Using eq.(2.1) one has
$$
\delta^2 F_{k+1, m} \left( \eta_{k+1}, \dots, \eta_m \right) = \ <-b_m
\eta_{m+1}, \eta_m>\ =\ <A^{(k)}_m \eta, \eta>,
$$
which contradicts the positivity of $\delta^2 F$.

We claim that the limit
\begin{equation}
\lim_{k \rightarrow - \infty} A^{(k)}_n = A_n
\label{24}
\end{equation}
exists and $A_n$ is a positive definite matrix sequence with the
recursion rule
\begin{equation}
A_{n+1} = a_{n+1} - b^T_n A^{-1}_n b_n .
\label{25}
\end{equation}
Indeed, it is easy to see by induction that $A^{(k)}_n$ is
monotone in $k$: $A_n^{(k)} > A_n^{(k-1)}$, for all $n>k$.  The
initial step $A^{(k)}_{k+1} > A^{(k-1)}_{k+1}$ follows from

$A^{(k)}_{k+1} = a_{k+1}$ and $A^{(k-1)}_{k+1} = a_{k+1} - b^T_k
\left[
  A^{(k-1)}_k \right]^{-1} b_k$, so $A^{(k)}_{k+1} - A^{(k-1)}_{k+1} = b^T_k
\left[ A^{(k-1)}_k \right]^{-1} b_k$.

The induction step is also simple: if
$$
A^{(k)}_n > A^{(k-1)}_n
$$
then
$$
A^{(k)}_{n+1} - A^{(k-1)}_{n+1} = - b^T_n \left( [A^{(k)}_n]^{-1} -
  [A^{(k-1)}_n]^{-1} \right) b_n.
$$

Thus the limit (2.4) exists and is a non-negative definite matrix.
Moreover $A_n$ is positive definite since it is necessarily
non-degenerate (together with $A^{(k)}_n$, the limit $A_n$ has to
satisfy the recurrence relation (2.5) which can be written without
the inverses of $A_n$). The claim is justified. We summarize the
result in the following

\begin{theorem}
\label{T4}
For any strict local minimal configuration $\left\{ q_n \right\}$ there
exists a
non-singular solution $\xi$ of the matrix Jacobi equation such that the
matrices $A_n = -b_n \xi_{n+1} \xi^{-1}_n$ are symmetric positive
definite and satisfy
\end{theorem}
\begin{equation}
A_{n+1} = a_{n+1} - b^T_n A^{-1}_n b_n . \label{26}
\end{equation}

\section{Proofs of the main theorems}

In this section we use the construction of the previous section to
prove Theorems {\ref{T1}} and {\ref{T2}}, and then apply them to prove
Theorem {\ref{T3}}.

\noindent{\bf Proof of Theorem {\ref{T1}}}

Consider the evolution transformations $R_m^n$ defined above and
the orbit of the point $x,x_n=R_0^nx$, and consider the
corresponding extremal $q_n=\pi x_n$.

Define ${\cal W}_n(x) = \lim_{k \rightarrow - \infty} {\cal
W}_n^{(k)} (x)$, where \quad $ {\cal W}_n^{(k)} (x) = (R_k^n)_*
({\cal V}(R_n^k x))$. Note, that by the assumption of no
conjugate points the Lagrangian subspaces ${\cal W}_n^{(k)} (x)$
are transversal to the vertical subspaces ${\cal V}(x)$. Moreover,
one can easily check that the corresponding matrices
$W_n^{(k)}$ satisfy:
\begin{equation}
W_n^{(k)} (x) = - \partial_{11} S_n (\pi(x), \pi(T_n(x)) +
A^{(k)}_n. \label{31}
\end{equation}

Therefore, by the properties of $A_n^{(k)}$ of the previous
section, the matrices $W_n$ are well defined and satisfy the equation:
\begin{equation}
\lim_{k \rightarrow - \infty} W_n^{(k)} (x) = W_n(x) = -
\partial_{11} S_n (\pi(x), \pi(T_n(x)) + A_n. \label{32}
\end{equation}
Notice that $W_n(x)$ depends measurably on $x$, since for every
$n,k$, ${\cal W}_n^{(k)}(x)$ is a smooth field of Lagrangian
subspaces.  The invariance property of the fields ${\cal W}_n$
follows immediately from the transformation rule

$${\cal W}_{n+1}^{(k)}(x)=(T_n)_* {\cal
W}_n^{(k)}(T_n^{-1}x)$$

\noindent for ${\cal W}_n^{(k)}$ which is immediate from the definition.
This yields the proof of theorem {\ref{T1}}.$\Box$

\noindent{\bf Proof of Theorem {\ref{T2}}} As in the proof of
Theorem 1 consider the orbit of the point $x$.
In order to prove
the inequalities 1 and 2 of Theorem 2, we shall use strongly that
all the matrices $A_n$ are positive definite. Then (3.1), (3.2)
imply
\begin{equation}
\nonumber -\partial_{11} S_n (\pi(x), \pi(T_n(x)) \leq
W_n(x).\label{33}
\end{equation}
And therefore
$$
-\partial_{11}S_n (q,q_{+}) \leq W_n(x).
$$
Also, using the relation (2.5), we have
\begin{eqnarray}
 \nonumber W_{n+1}(T_nx) = -\partial_{11} S_{n+1} (\pi (T_nx), \pi
(T_{n+1}\circ T_nx)) + A_{n+1}=  \\
 =-\partial_{11}S_{n+1} (\pi
(T_nx), \pi
(T_{n+1}\circ T_nx))  + a_{n+1} - b_n^T A_n^{-1}b_n =\\
\label{34}
\nonumber
\partial_{22} S_n (\pi x, \pi (T_n x)) - b_n^T A_n^{-1}b_n
 \leq \partial_{22} S_n (\pi x, \pi (T_n x)).
\end{eqnarray}

Thus we have
\begin{equation}
W_n(x) \leq \partial_{22} S_{n-1} (\pi (T^{-1}_{n-1}x), \pi (x))=
\partial_{22} S_{n-1}(q_{\_}, q).
\label{35}
\end{equation}
Notice that the inequalities (3.3) and (3.5) can be expressed
geometrically by
$$(T_{n+1}^{-1})_* {\cal V}(T_nx) \le {\cal W}_n(x) \le (T_{n-1})_*
\left({\cal V}(T_n^{-1}x) \right)$$
This proves the first part of Theorem 2.

In order to prove the second part we subtract the two expressions
(3.4) and (3.2) for $W.$   We have
\begin{equation}
W_{n+1} (T_n x) - W_n (x) = \partial_{22} S_n (\pi x, \pi (T_n
x))+\partial_{11} S_n (\pi x, \pi (T_n x))-A_n- b_n^T A_n^{-1}b_n
\label{36}
\end{equation}
This can be rewritten as
\begin{eqnarray}
\nonumber
 W_{n+1} (T_n x) - W_n (x)= \partial_{11} S_n (q, q_{+})
+ \partial_{22} S_n (q, q_{+})-\\
- \left( A^{1 \over 2}_n + b^T_n A^{- {1\over 2}}_n \right) \left(
A^{1 \over 2}_n + A^{- {1\over 2}}_n b_n \right) + b_n + b^T_n .
\label{37}
\end{eqnarray}
Notice that the first matrix in brackets of (3.7) is the transpose
of the second one and thus
\begin{eqnarray}
\nonumber W_{n+1} (T_n x) - W_n (x) \leq \partial_{11} S_n (q,
q_{+}) +\partial_{22} S_n (q, q_{+})+b_n+b_n^T=\\
= \partial_{11} S_n (q, q_{+}) + \partial_{22} S_n (q, q_{+})
+ \partial_{12} S_n (q, q_{+}) + \partial_{21} S_n (q, q_{+}).
\label{38}
\end{eqnarray}
Moreover the inequality (3.8) is strict except when
\begin{equation}
A_n = - b_n = - b^T_n. \label{39}
\end{equation}
In the last case the expressions for $W_{n+1} (T_n x)$ and $W_n
(x)$ are
\begin{eqnarray}
W_{n+1} (T_n x)&=& \partial_{22} S_n (q, q_{+}) + \partial_{12} S_n (q, q_{+})
\nonumber \\
W_n(x) &=& - \partial_{11} S_n (q, q_{+}) - \partial_{21} S (q,
q_{+}) . \label{310}
\end{eqnarray}
This finishes the proof of Theorem 2. $\Box$

\noindent{\bf Proof of Theorem {\ref{T3}}}

In the case of a Frenkel-Kontorova functional we have
$$S_n(q,Q)=\frac{1}{2} (Q-q)^2+V_n(q),$$ where $V_n$ is periodic in
$q.$ In this case the partial derivatives of $S_n$ are
\begin{eqnarray}
  \nonumber \partial_{22} S_n=-\partial_{12} S_n=-\partial_{21} S_n=I\\
  \partial_{11} S_n=I+Hess(V_n) .
\end{eqnarray}
Suppose that all the extremals of the Frenkel-Kontorova functional
are without conjugate points. Then construct the fields of
Lagrangian subspaces ${\cal W}_n$ and the corresponding matrix
functions $W_n$ as in Theorems 1,2. Let us define
$$w_n(x)=tr W_n(x)$$
then $w_n$ is a bounded measurable function satisfying the
following inequality (a consequence of Theorem 2)
$$
w_{n+1}(T_nx) - w_n(x) \leq tr \left( \partial_{11} S_n (q, q_{+})
+ \partial_{22}
 S_n (q, q_{+}) + 2 \partial_{12} S_n (q, q_{+})  \right).$$
In other words we get the following
\begin{equation}
w_{n+1}(T_nx) - w_n(x) \leq \Delta V_n(q) .
\end{equation}
We shall see below that if all the extremals of the Frenkel-Kontorova
functional have no conjugate points then for almost
all $x$ there is equality in (3.12). Therefore by
Theorem 2, (3.10) holds, i.e. by the formulae (3.11)

$$W_n=-HessV_n \mbox{ and }  W_{n+1}(T_nx)=0.$$
In other words
$$W_n\equiv -Hess(V_n) \equiv 0$$ for all $n$.  But then
all the functions $V_n$ are constant.  This will finish the
proof of Theorem 3.

In order to establish equality in (3.12) we shall consider two
cases. In the first case the sequence $V_n$ is periodic, i.e.
$V_{n+p}\equiv V_n$ for some positive integer $p$ and for all $n.$
In the second case the sequence $V_n$ is of compact support, i.e.
$V_n\equiv const$ for $|n|>N$ for some $N.$ Consider first the
periodic case. In this case obviously ${\cal W}_{n+p} \equiv {\cal
W}_n$ and thus $w_n \equiv w_{n+p}$. Now we apply (3.12) $p$ times
to obtain
\begin{eqnarray}
\nonumber w_{n+p}(T_{n+p-1} \circ \cdots \circ
T_nx)-w_n(x)\leq\Delta V_n(\pi x)+\\+\Delta V_{n+1}(\pi (T_n
x))+\ldots+\Delta V_{n+p-1}( \pi(T_{n+p-2} \circ \cdots \circ T_n
x)) .
\end{eqnarray}

Let us recall the additional property of the standard maps $T_n$
that the phase space is effectively compact (see remark in Example
1). This implies immediately that  each field ${\cal W}_n(x)={\cal
W}_n (p,q)$ depends periodically on $p$ as well as on $q$.  Thus
the function $w_n$ is a periodic bounded function on $\Pi^{2d}$.
Now we can finish the argument by the following reasoning. If
there is strict inequality in (3.12) for some $n$ on a set of
positive measure then one has strict inequality in (3.13) also on
a set of positive measure. But then the strict inequality remains
after the integration of (3.13) over the whole phase space $ \Pi
^{2d}$. But this is a contradiction, because since all the
transformations $T_n$ are symplectic (and hence measure
preserving) then one can easily see that the integrals of both
sides of (3.13) over $ \Pi ^{2d}$ vanish. This finishes the proof
of the claim in the periodic case.

In the second case the idea is similar. The important ingredient
in its realization is the following claim. The limit
$$\lim_{n \rightarrow \pm \infty} w_n(x)=0$$ exists and the convergence
is uniform in $x$. In fact for those $n$ which lie to the left of
the support of $V_n$ it easily follows from the construction that
$W_n(x)=0$ and then $w_n(x)=0$ for all $x$. For large positive $n$
we have $V_n \equiv 0$, thus the recursion rule for the matrices
$W_n,A_n$ from (2.6) and (3.2) is:
\begin{equation}
A_{n+1}=2I-A_n^{-1} \mbox{ and } W_n=-I+A_n .
\end{equation}
Then for the eigenvalues of $A_n$ the same recursion rule holds $$
\lambda_{n+1}=2-\frac{1}{\lambda _n}.$$ Recall that all the
matrices $A_n$ are positive definite. Therefore all $\lambda_n$ are
positive and then one can easily see that the sequence $\lambda
_n$ is monotonically decreasing and converges to $1$. Moreover,
it is clear from the formula that $\lambda _{n+1}$ is less than
$2$. Therefore, $A_n$ converges (uniformly for all
orbits) to I and thus $W_n$ to $0.$ This proves the claim. In
order to finish the proof of the Theorem
one proceeds exactly as in the previous
case. One takes $N$ sufficiently large and sums up the inequality
(3.12) from $-N$ to $N$. This completes the proof of theorem 4.
 $\Box$

\section{Discussion and some open questions}

\begin{enumerate}
\item The variational principle (1.3) can be considered on other
  configuration manifolds different from tori, for example on
  hyperbolic manifolds.  It would be interesting to understand the
  consequences of the no conjugate points condition
  for these cases.  Another very interesting direction would be to
  study, along the lines of this paper, variational principles of the
  form (1.3) on configurations $\left\{ q_n \right\}$ for $n$ lying on
  some lattice $\Bbb Z^k$ (see also \cite{[K-L-R]}).  Some results in
  this direction were obtained in \cite{[B-M]} for
  multi-continuous-time systems.
\item An important problem is to understand to what extent
the smoothness of $W$ is required.
  An example of not smooth enough $W$ would give a qualitatively new
  system without conjugate points.
\item The integration trick used in the proof of Theorem {\ref{T3}}
  worked well due to compactness of the phase space for the standard
  map.  In many interesting cases, however, the phase space is not
  compact.  Then new integral-geometric approaches are required.  For
  example it is not clear yet how to apply this to the so-called outer
  billiard problem \cite{[Ta]}.  It would be reasonable to
  conjecture that the only outer billiards without conjugate points on
  the affine plane are the elliptic ones.  In some
  cases the lack of compactness can be overcome \cite{[B-P1],[B2]}.
\item It was proved by J Moser \cite{[M]} for area-preserving twist
  maps that every such map can be seen as the time-one map of an
  optical Hamiltonian function.  This result was generalized in
  \cite{[B-P2]} to higher dimensions for those twist maps with
  symmetric matrix $\partial_{12} S$ (see \cite{[Go]} for the proof
  and discussion).  It is not clear what can be said about the
  interpolation problem for symplectic twist maps without conjugate
  points.  Is it true that they can be interpolated by flows without
  conjugate points?
\item One might prefer an extended notion of conjugate points for symplectic
twist maps, which mimics more closely the properties of maps arising from
optical Hamiltonian flows, by allowing a conjugate point to occur in between
two integer times.  To formalise this, we say that an orbit of Lagrange
planes {\em crosses the vertical} between times $n$ and $n+1$ if the
signature of the associated quadratic form changes.  Then for $m<n$
we can say time
$m$ is conjugate to $(n,n+1)$ along orbit $(x_i)$ if the orbit of the
plane which is vertical at time $m$ crosses the vertical between times
$n$ and $n+1$.  Similarly for $m > n+1$ by using the backwards dynamics.
Also we can say $(m,m+1)$ is conjugate to $(n,n+1)$ if the orbit of the
vertical plane at time $m$ crosses the vertical between times $m,m+1$ and
between times $n,n+1$.  The definition of this paper is incorporated by
saying times $m$ and $n$ are conjugate if the orbit of the vertical at
time $m$ has non-zero intersection with the vertical at time $n$.  Of
course, if all orbits have no conjugate points in this extended sense then
they have no conjugate points in the restricted sense and hence the
conclusions of the paper still follow.  Advantages of the extended
definition are that possession of conjugate points becomes stable and that
for discretisations of an orbit of an optical Hamiltonian system its conjugate
points are inherited.
\end{enumerate}

\hfill


\begin{thebibliography}{References}
\bibitem{[B1]} Bialy, M., {\em Convex Billiards and a theorem by
  E. Hopf}, Math. Z. {\bf 24} (1993) 147--154.

\bibitem{[B2]} Bialy, M., {\em On shocks formation in forced Burgers
  equation and application to a quasi-linear system}, GAFA, Vol.10
(2000) 732-740.

\bibitem{[B-M]} Bialy, M., MacKay, R.S., {\em Variational properties of
  a non-linear elliptic equation and rigidity}.
  Duke Math. J. 102 (2000) 391--401.

\bibitem{[B-P1]} Bialy, M., Polterovich, L., {\em Hopf type rigidity
  for Newton equations}, Math. Research Letters {\bf 2} (1995)
  695--700.

\bibitem{[B-P2]} Bialy, M., Polterovich, L., {\em Hamiltonian systems,
  Lagrangian tori and Birkhoff theorem}, Math. Ann. {\bf 292} (1992)
  619--627.

\bibitem{[C-F]} Croke, C., Fathi, A., {\em An inequality between energy
    and intersection}, Bull. Lond. Math. Soc. {\bf 22}, (1990)
    489--494.

\bibitem{[C-I]} Contreras, G., Iturriaga, R., {\em Convex Hamiltonians
    without conjugate points}, Ergodic Theory Dynam. Systems 19 (1999), 
no. 4, 901--952

\bibitem{[Go]} Gol\'e, C., {\em Symplectic twist maps}, Advanced
    series in Nonlinear Dynamics, Vol. 18, (World Scientific, 2001).

\bibitem{[G]} Green, L., {\em A theorem of E. Hopf},
  Mich. Math. Journal {\bf 5}, (1958) 31--34.

\bibitem{[He]} Herman, M., {\em Inegalites a priori pour des tores
    lagrangiens invariants par des diff\'eomorphismes symplectiques},
    Pub. Math. IHES {\bf 70}, (1990).

\bibitem{[K]} Knauf, A., {\em Closed orbits and converse KAM theory},
  Nonlinearity {\bf 3}, (1990) 961--973.

\bibitem{[K-L-R]} Koch, H., de la Llave, R., Radin, C., {\em
  Aubry-Mather theory for functions on lattices}, Discrete
  Contin. Dynam. Systems {\bf 3}, (1997) 135--151.

\bibitem{[M-M-S]} MacKay RS, Meiss JD, Stark J, {\em Converse KAM
theory for symplectic twist maps}, Nonlinearity 2 (1989) 555--570.

\bibitem{[M]} Moser, J., {\em Monotone twist mappings and the calculus
  of variations}, Erg. Th. and Dyn. Sys. {\bf 6}, (1986) 401--413.

\bibitem{[M-V]} Moser, J., Veselov, A., {\em Discrete versions of some
  classical integrable systems and factorization of matrix
  polynomials}, Comm. Math. Physics {\bf 139}, (1991), 217--243.

\bibitem{[Ta]} Tabachnikov, S., {\em On the dual billiard problem},
  Adv. Math. {\bf 115}, (1995), 221--249.

\bibitem{[W]} Wojtkowski, M., {\em Two application of Jacobi fields to
  the billiard ball problem}, J. Diff. Geom. {\bf 40}, (1994),
  155--164.





\end{thebibliography}
\end{document}